\documentclass[%review,
12pt,
authoryear,
times
]{elsarticle}

\usepackage{amsmath}

\usepackage[colorlinks=true, %pdfstartview=FitV, linkcolor=black, citecolor= black, urlcolor= blue
]{hyperref}
\usepackage{lineno}

\journal{Advances in Space Research}

\begin{document}

\begin{frontmatter}

\title{ LEO intermediary propagation as a feasible alternative to Brouwer's gravity solution\tnoteref{KePASSA} }
\tnotetext[KePASSA]{Presented in KePASSA 2014, Logro\~no, Spain, April 23--25, 2014}

\author%[roa]
{ Martin Lara\fnref{GRUCACI} %\corref{cor1} 
}
\fntext[GRUCACI]{ GRUCACI, University of La Rioja, Spain }

\ead{mlara0@gmail.com}

\address%[roa]
{Columnas de H\'ercules 1, ES-11000 San Fernando, Spain}


\date{}

%\maketitle

\begin{abstract}
The performance of Brouwer's gravity solution is compared with Deprit's second-order radial intermediary. Taking the main problem of the artificial satellite as the test model, that is perturbations are limited to the effects of the second zonal harmonic, it is demonstrated that the intermediary solution provides an efficient alternative for the analytical propagation of low earth orbits in a range determined by non-impact orbits with eccentricities below one tenth.
\end{abstract}


\end{frontmatter}

%\linenumbers

\section{Introduction}

Increasing needs of collision avoidance stir recent concern in reviewing SGP4 \citep{HootsRoehrich1980}, an orbit prediction model that is customarily used in the propagation of two-line element sets \citep[see also][and references therein]{HootsSchumacherGlover2004,ValladoCrawfordHujsakKelso2006}. %\footnote{}.
Among the reasons for this concern is the need of including uncertainties estimation in the predictions, but also the detected significant along-track errors which may be related to missing terms in the SGP4 implementation of the gravitational theory \citep{Kelso2007,Easthope2014}. These facts motivated me to search the literature for alternatives to the gravitational solution used by SGP4. I found that these alternatives do exist, can be refined to cope with more subtle short-period effects, and may provide increased efficiency in terms of precision and computation time when compared with the classical approach. In spite of comparisons are limited to the main problem of artificial satellite theory, the proposed alternative may be further explored to include more disturbing effects which are required for realistic propagation of space debris or orbits about asteroids.
%
%In particular I reimplemented Brouwer's solution for the main problem of artificial satellite theory and extended Deprit's radial intermediary by computing second order periodic corrections.}
\par

Brouwer's celebrated closed form solution of the earth's artificial satellite problem \citep{Brouwer1959} is in the roots of SGP4. 
In spite of the undeniable merits of Brouwer's seminal approach, or its variants for properly dealing  with small eccentricities and inclinations \citep{Lyddane1963,CohenLyddane1981}, proceeding by averaging is not the unique possibility in obtaining analytical solutions by general perturbations, and the use of \emph{intermediary} orbits like useful approximate solutions of the problem of artificial satellite theory (AST) was proposed as early as in 1957 \citep{Sterne1958}. 
\par

AST intermediaries are commonly obtained by reorganizing the terms of the disturbing function, a simple expedient that may be preceded by the simultaneous addition and subtraction to the geopotential of some smartly-chosen supplementary terms. After reorganization, a part of the Hamiltonian that admits a separable generating function is taken as the (zero order) integrable problem whereas the rest of the disturbing function is taken as the perturbation, which may be further neglected \citep[see][for instance]{Aksnes1965,Oberti2005}. Hence, \emph{common} intermediaries are formulated in the same variables as the original satellite problem ---traditionally in spherical variables.
\par

It is habitually accepted that useful intermediary solutions should retain all the first-order secular and long-period effects of the artificial satellite problem and as much as possible of the short-period effects \citep[cf.][]{GarfinkelAksnes1970}. Because of that, neither the Keplerian Hamiltonian nor the equatorial main problem\footnote{We call the \emph{equatorial main problem} to the perturbed Keplerian motion of a satellite in the main problem assumptions that is constrained to lie on the earth's equatorial plane.} are considered suitable intermediaries in spite of their integrability \citep{Jezewski1983}. Including some of the second order secular effects of the geopotential in the intermediary is highly desirable too, and hence solutions devised by \citet{Vinti1959,Vinti1961,Vinti1966} or equivalent ones based on the generalized problem of two fixed centers \citep[see][and references therein]{AksenovGrebenikovDemin1962,LukyanovEmeljanovShirmin2005} are considered remarkable achievements.\footnote{Vinti found a particular choice for the coefficients of the zonal potential
% $J_{2m}=(-1)^{m+1}J_2^m$ 
that makes it separable in oblate spheroidal coordinates. Vinti's analytical solution closely applies to earth's orbits where the zonal coefficient of degree 4 almost matches his choice. The equivalence between Vinti's model and the problem of two fixed centers was later recognized.}
\par

Eventually, the efforts of \citet{CidLahulla1969} in obtaining a competitive alternative to Brouwer's solution produced a major breakthrough in the search for efficient intermediaries. Proceeding in polar-nodal variables, the usual polar ones complemented with the argument of the ascending node and its conjugate momentum, Cid and Lahulla showed that, contrary to the two canonical transformations required in Brouwer's approach, a single contact transformation is enough for removing the argument of the latitude from the main problem Hamiltonian, which after that, and in view of the cyclic character of the argument of the node, results to be integrable.
\par

Up to the first order, Cid and Lahulla's Hamiltonian is formally equal to the radial ---and, therefore, integrable--- part of the main problem Hamiltonian, in this way highlighting its relationship with common intermediaries, on one hand, and disclosing the important role played by polar-nodal variables in the search for them, on the other. Besides, it incorporates all the first-order secular effects of the main problem \citep{DepritFerrer1987} as well as the short-period terms that affect the radial distance, thus fulfilling the traditional requisites for acceptable main problem intermediaries. Furthermore, Cid and Lahulla's solution also takes into account the remaining first-order periodic terms of the main problem, which are recovered by means of a contact transformation. All these facts inspired Deprit, who introduced the concept of \emph{natural} intermediaries, which are integrable after a contact transformation, and demonstrated that most of the common intermediaries can be ``naturalized'' by finding the transformation that turns the main problem into the intermediary \citep{Deprit1981}. Remarkably, Deprit also showed how these contact transformations may be extended to higher orders, and proposed his own radial intermediary which, as opposite to other existing intermediaries, does not rely on the evaluation of elliptic functions. Later efforts based on Deprit's approach, showed how the generation of intermediaries of AST can be systematized \citep{FerrandizFloria1991}.
\par

In spite of the efforts in improving the intermediaries' performance by including second order effects of the gravity potential \citep{Garfinkel1959,Aksnes1967,AlfriendDasenbrockPickardDeprit1977,DepritJGCD1981,CidFerrerSeinEchaluce1986,DepritFerrer1989}, increasing requirements on satellite orbit prediction accuracy soon led to a decline in interest in intermediaries, to favor instead analytical and semi-analytical theories based on averaging. Still, the use of intermediary orbits of AST is enjoying a revival these days, and new applications of intermediary-based solutions had been recently encouraged for different purposes, as for predicting long-term bounded relative motion \citep{LaraGurfil2012} or like a suitable choice for onboard orbit propagation as opposite to the usual numerical integration \citep{GurfilLara2014}.
\par

Here, the main problem of AST is taken as a test model to demonstrate that the use of intermediary orbits may provide an efficient alternative to Brouwer's gravity solution for the propagation of low earth orbits (LEO) in the short time intervals required by usual catalogue maintenance. Indeed, while natural intermediary closed-form solutions are unquestionably of higher complexity than Brouwer's secular terms and, in consequence, their evaluation becomes much more computationally costly, they only rely on one simplification transformation whose terms are definitely simpler than the corresponding Fourier series required by Brouwer's double averaging approach. Besides, AST intermediaries do not suffer from mathematical singularities at the critical inclination \citep[see reviews in][]{CoffeyDepritMiller1986,Jupp1988,Lara2014Rome} and hence do not need to rely upon the functional patches on which higher-order analytical theories by averaging unavoidably depend ---cf.~Section 7 of \citep{CoffeyNealSegermanTravisano1995} or appendix A.F.~in \citep{HootsSchumacherGlover2004}.\footnote{It must be emphasized that neither analytical theories based on a double averaging nor existing intermediaries can deal properly with second-order perigee effects, which are crucial in the study of frozen-perigee orbits. This special class of orbits is usually approached either from a semi-analytical perspective or a purely numerical one \citep{CuttingFrautnickBorn1978,Broucke1994,CoffeyDepritDeprit1994,Lara1999,Lara2008}.}
\par

On the other hand, efforts in solving higher-order intermediaries by separation of variables have been unsuccessful, in spite of their integrable character. This shortcoming imposes general perturbations algorithms to progress with an additional transformation in order to achieve the analytical solution of the intermediary to higher orders \citep{Aksnes1967,DepritJGCD1981,DepritRichardson1982}. However, this extra transformation may well be avoided in the case of LEO, an instance in which orbital eccentricity is small. % thus pushing the perigee dynamics to higher orders. 
In particular, this is the case of Deprit's radial intermediary, or DRI in short \citep{Deprit1981}. For the lower-eccentricity orbits, which are a vast majority in a space catalogue of earth orbiting objects, DRI  can be solved up to higher-order effects by separation of the generating function \citep{DepritJGCD1981,Floria1993}. 
%More precisely, the intermediary is reduced to a quasi-Keplerian system which retains higher-order effects of the main problem solution and admits standard closed-form solution.
\par

In the present research, a quasi-Keplerian intermediary has been constructed which admits standard closed-form solution and is of straightforward evaluation. The proposed intermediary is based on a simplification of DRI which is valid for low-eccentricity orbits and has been improved with the computation of those second order short-period corrections that are consistent with the simplification assumptions of the theory, in this way incorporating second-order secular as well as periodic terms of the main problem into the analytical solution. Also, a complete reordering of the periodic corrections has been carried out that notably fastens their evaluation. 
\par

Taking the numerical integration of the main problem as the true solution, a variety of tests have been carried out in a representative set of LEO with inclinations encompassing from equatorial to polar. These tests show that the performance of the second-order intermediary is clearly better than that of Brouwer's analogous propagations in nonsingular variables for short time intervals. Besides, the intermediary-based LEO propagator remains very competitive for time intervals spanning up to several weeks.
\par

Propagation of space debris and orbits about asteroids undoubtedly requires more complete forces model than the one discussed here. In particular orbit propagation about highly inhomogeneous gravity fields definitely requires inclusion of tesserals effects, which can also be dealt with analytically \citep[see][and references therein]{Garfinkel1965,CeccaroniBiggs2013,LaraSanJuanLopezOchoa2013a} albeit the efficiency of the closed form approach may be questionable in some cases \citep[cf.][]{LaraSanJuanLopezOchoa2013}. The usefulness of replacing parts of the standard perturbation solution by analogous results based on intermediary orbits still remain to be explored.

\section{Deprit's radial intermediary}

Assume an inertial reference system defined by the earth's center of mass, the $z$ axis defined by the earth's rotation axis, and the $x$ and $y$ axes lying in the equatorial plane and defining a direct frame.
In the canonical set of polar-nodal variables $(r,\theta,\nu,R,\Theta,N)$ ---standing for the radial distance from the origin, the argument of the latitude, the longitude of the ascending node, the radial velocity, the modulus of the angular momentum vector, and the projection of the angular momentum over the $z$ axis, respectively--- the main problem Hamiltonian is written
\begin{equation} \label{hmp}
\mathcal{H}=\frac{1}{2}\left(R^2+\frac{\Theta^2}{r^2}\right)-\frac{\mu}{r}\left[
1-J_2\,\frac{\alpha^2}{r^2}\,P_2(s\sin\theta)\right],
\end{equation}
where $\mu$ is the earth's gravitational parameter, the scaling factor $\alpha$ is the earth's mean equatorial radius, and $J_2$ is the second zonal harmonic coefficient of the geopotential. Besides, $P_2$ is the Legendre polynomial of degree 2, $s\equiv\sin{i}$, and $c\equiv\cos{i}=N/\Theta$. In view of $\nu$ is ignorable in Eq.~(\ref{hmp}), its conjugate momentum $N$ is an integral of the motion, reflecting the symmetry of the main problem dynamics with respect to rotations about the $z$ axis.
\par

Deprit's radial intermediary is obtained after simplifying the main problem Hamiltonian by applying to it the elimination of the parallax, which is a canonical transformation 
\begin{equation*}
(r,\theta,\nu,R,\Theta,N)\longrightarrow(r',\theta',\nu',R',\Theta',N')
\end{equation*}
that removes non essential short-periodic terms of the original Hamiltonian without reducing the number of degrees of freedom \citep{Deprit1981}.
\par

The elimination of the parallax from the main problem Hamiltonian is fully documented in the literature so details are not provided here. The interested reader may consult the original paper of \citet{Deprit1981} or my simpler re-derivation of this canonical simplification based on the use of Delaunay variables \citep{LaraSanJuanLopezOchoa2013b}.
\par

After eliminating the parallax, and up to the second order of $J_2$, the main problem Hamiltonian in new (prime) polar-nodal variables is %\citep[cf.][]{Deprit1981}
\begin{equation} \label{hpnp}
\mathcal{H}=H_0+H_1+\frac{1}{2}H_2+\mathcal{O}(J_2^3)
\end{equation}
where, dropping primes for brevity,
\begin{align} \allowdisplaybreaks \label{H00}
H_0 \;=&\; \frac{1}{2}\left(R^2+\frac{\Theta^2}{r^2}\right)-\frac{\mu}{r}, \displaybreak[0] \\ \label{H01}
H_1 \;=&\;-\varepsilon\,\frac{\Theta^2}{r^2}\left(1-3c^2\right) \displaybreak[0] \\ \label{H02}
H_2 \;=&\; \varepsilon^2\,\frac{\Theta^2}{r^2}\left\{ 1-21c^4-6\left(c^2-\frac{5}{8}s^4\right)(\kappa^2+\sigma^2) \right. \\ \nonumber
& \left. +21\left(1-\frac{15}{14}s^2\right)s^2\left[(\kappa^2-\sigma^2)\cos2\theta+2\kappa\,\sigma\sin2\theta\right] \right\}.
\end{align}
For abbreviation, it has been introduced the notation
\begin{equation}
\varepsilon\equiv\varepsilon(\Theta)=-\frac{1}{4}J_2\frac{\alpha^2}{p^2},
\end{equation}
where $p$ is the parameter of the conic, or \textit{semilatus rectum},
\begin{equation}
p=\Theta^2/\mu,
\end{equation}
and $\kappa$ and $\sigma$ are the projections of the eccentricity vector in the orbital frame % radial and normal directions
\begin{equation} \label{ksef}
\kappa=e\cos{f}, \qquad \sigma=e\sin{f},
\end{equation}
where $f$ is the true anomaly and $e$ is the eccentricity. Based on the conic solution of the Keplerian motion, these functions are easily expressed in polar-nodal variables as
\begin{equation} \label{ef}
\kappa\equiv\kappa(r,\Theta)=\frac{p}{r}-1, \qquad \sigma\equiv\sigma(R,\Theta)=\frac{p\,R}{\Theta}.
\end{equation}
Note that $f\equiv{f}(r,R,\Theta)$ and $e\equiv(r,R,\Theta)$.
\par

The appearance of $\theta$ in Eq.~(\ref{H02}) prevents integrability, and the transformed Hamiltonian in Eq.~(\ref{hpnp}) has the same degrees of freedom as the original main problem in addition to being yet more intricate. However, the new Hamiltonian is said to be ``simplified'' \citep{DepritFerrer1989}. Indeed, because the argument of perigee is $\omega=\theta-f$, it is simple to derive from Eq.~(\ref{ksef}) the relations
\begin{eqnarray} \label{ecg}
e\cos\omega %=e\cos(f-\theta) 
&=& \kappa\cos\theta+\sigma\sin\theta, \\ \label{esg}
e\sin\omega %=e\sin(f-\theta) 
&=& \kappa\sin\theta-\sigma\cos\theta,
\end{eqnarray}
from which
\begin{eqnarray} 
e^2\cos2\omega &=& (\kappa^2-\sigma^2)\cos2\theta+2\kappa\sin2\theta,
\\ \label{e2}
e^2 &=& \kappa^2+\sigma^2.
\end{eqnarray}
Hence, Eq.~(\ref{H02}) is rewritten as
\begin{equation} \label{H02d}
H_2=\varepsilon^2\,\frac{\mu}{r}\,\frac{p}{r}\left[1-21c^4-6\left(c^2-\frac{5}{8}s^4\right)e^2
+21\left(1-\frac{15}{14}s^2\right)s^2e^2\cos2\omega\right],
\end{equation}
in this way showing that the maximum negative power of $r$ in Eq.~(\ref{hpnp}) is $-2$ contrary to the power $-3$ in Eq.~(\ref{hmp}). This fact makes a subsequent closed form elimination of the short-period terms from Eq.~(\ref{hpnp}) by means of the usual Delaunay normalization \citep{Deprit1982} much simpler than Brouwer's direct elimination, on the one hand, and eases the construction of higher order solutions, on the other.\footnote{The total number of terms of the short-period transformations resulting from the elimination of the parallax and a consequent normalization is claimed to be reduced to just one fourth of the number of terms required in a third-order solution in closed when computed by the classical, Brouwer style normalization \citep[cf.][]{CoffeyDeprit1982}.}
\par

In spite of the non-integrable character of the simplified Hamiltonian, if terms multiplied by $e^2\,J_2^2$, which will be of higher order when $e^2\sim{J}_2$ (say $e<0.1$ for the earth if we neglect terms of the order of $10^{-8}$), are neglected then Eq.~(\ref{H02}) is converted into
\[
H_2=\frac{\Theta^2}{r^2}\,\varepsilon^2\left(1-21c^4\right),
\]
and the argument of the latitude becomes cyclic as it is cancelled from Eq.~(\ref{H02d}). Therefore, Eq.~(\ref{hpnp}) is simplified to the radial intermediary
\begin{equation} \label{qks}
\mathcal{D}=\frac{1}{2}\left(R^2+\frac{\tilde\Theta^2}{r^2}\right)-\frac{\mu}{r},
\end{equation}
where
\begin{equation}
\tilde\Theta=\Theta\,\sqrt{1-\varepsilon\,(2-6c^2)+\varepsilon^2\,(1-21c^4)}.
%=\Theta\left[1-\varepsilon\,(1-3c^2)+3\varepsilon^2(1-5c^2)c^2\right]
\end{equation}
Equation (\ref{qks}) is a \emph{quasi-Keplerian} system with varying (depending on $J_2$) angular momentum $\tilde\Theta$, which is called here DRI.\footnote{In fact, in view of Eqs.~(\ref{e2}) and (\ref{ef}), the eccentricity does not depend on $\theta$ and integrability can be obtained by neglecting only from Eq.~(\ref{H02d}) the terms that are affected by $\omega$, in this way keeping in the intermediary more effects of the original problem. But the integration of this more complete intermediary requires the apparatus of perturbation theory \citep[cf.][]{DepritJGCD1981} which is avoided in the case of quasi-Keplerian systems.}
\par

The integration of Eq.~(\ref{qks}) can be done by the standard Hamilton-Jacobi reduction to Delaunay-similar variables. This procedure requires the introduction of auxiliary variables for performing the quadratures, the so-called \emph{anomalies}, and solving the Kepler equation. For completeness, the necessary formulae are detailed in \ref{s:qK}.
\par

However, the quasi-Keplerian system provides only the solution in the prime space, and the integration of DRI is not completed until the short-period terms removed by the elimination of the parallax are recovered. Sequences for recovering these terms from the direct transformation, as well as the inverse transformation required for computing initial conditions in the prime space, are provided in \ref{s:shp}. Note that, in spite of the transformation from and to prime variables is computed in closed form of the eccentricity, in consequence with the assumption made of neglecting terms of $\mathcal{O}(e^2)$ from the second order Hamiltonian after the elimination of the parallax, the transformation equations for the second order short-period effects have been analogously simplified. Note also that evaluation of the corrections only involves sine and cosine functions of arguments $2\theta$ and $4\theta$.

\section{Comparison with Brouwer's solution}

Brouwer's gravity solution is accomplished by finding the canonical transformations that reduce the satellite problem to its secular terms \citep{Brouwer1959}. These transformations are carried out in the set of Delaunay variables $(\ell,g,h,L,G,H)$,
% which is the canonical counterpart of the traditional set of Keplerian elements $(a,e,i,\Omega,\omega,M)$, standing for semi-major axis, eccentricity, inclination, Right Ascension of the ascending node, argument of the periapsis, and mean anomaly, respectively. The Delaunay set 
which comprises three angles: the mean anomaly $\ell$, the argument of the perigee $g$, and the argument of the node $h$, as well as their  three canonical conjugate momenta: the Delaunay action $L=\sqrt{\mu\,a}$, conjugate to $\ell$, the total angular momentum $G=L\sqrt{1-e^2}$, conjugate to $g$, and the polar component of the angular momentum $H=G\cos{i}$, conjugate to $h$. Note that $g$ was previously noted $\omega$, and that $h$, $G$, and $H$ are noted $\nu$, $\Theta$, and $N$, respectively, in polar-nodal variables.
\par

The computation of ephemeris based on Brouwer's canonical transformations in Delaunay variables is known to introduce errors of the first order in the short-period corrections when the eccentricity approaches to zero. This happens because of the singularity of Delaunay variables in the case of zero eccentricity orbits.\footnote{Delaunay variables are also singular for zero inclination orbits} However, since these singularities are just virtual \citep{Henrard1974}, the trouble is avoided by expressing the perturbations of the elements in a different set of nonsingular variables.
\par

Thus, in order to make Brouwer's solution useful also for low-eccentricity and low-inclination orbits, instead of using Delaunay variables \cite{Lyddane1963} resorts to the set of Poincar\'e canonical elements
\begin{align} \allowdisplaybreaks
x_1 = & L \displaybreak[0] \\
x_2 = & \sqrt{2L}\,\sqrt{1-\eta}\cos(g+h) \displaybreak[0] \\
x_3 = & \sqrt{2L}\,\sqrt{\eta\,(1-c)}\cos{h} \displaybreak[0] \\
y_1 = & \ell+g+h \displaybreak[0] \\
y_2 = &-\sqrt{2L}\,\sqrt{1-\eta}\sin(g+h) \displaybreak[0] \\
y_3 = &-\sqrt{2L}\,\sqrt{\eta\,(1-c)}\sin{h} 
\end{align}
where
\begin{equation}
\eta=\sqrt{1-e^2}=\frac{G}{L}
\end{equation}
is commonly called the ``eccentricity function''.
\par

Note that, in Lyddane's view, the perturbation theory is not constructed in Poincar\'e elements, but in Delaunay ones. After that, the perturbation theory is reformulated in the desired set of variables, which, furthermore, do not need to be canonical. Indeed, other different sets of nonsingular variables can be used, and, because of the axial symmetry of the main problem, the evaluation of Brouwer's analytical solution is found more expedite when using the non-canonical, nonsingular set:
\begin{equation} \label{FCShHL}
F=\ell+g, \quad C=e\cos{g}, \quad S=e\sin{g}, \quad h, \quad H, \quad L,
\end{equation}
proposed by \cite{DepritRom1970}.
\par

The locution ``Brouwer's solution'' is loosely applied to different versions of the perturbation solution to the problem of the artificial satellite found by \citet{Brouwer1959}. Here, this term is used to refer to an analytical solution obtained by double averaging, which  comprises the secular terms up to the second order, but is limited to first-order periodic terms ---in agreement with common implementations of this theory. In particular, because the concern is limited to the propagation of low-eccentricity orbits, Brouwer's solution is formulated in the nonsingular variables proposed by \citet{DepritRom1970}. Besides, Brouwer's solution is constrained to the main problem, and since it is compared with DRI both in accuracy and computing time, in order to obtain results that are as far as possible unbiased Brouwer's solution has been simplified to match the assumptions made in the analytical integration of DRI. Namely, it has been truncated by neglecting the order of $\mathcal{O}(J_2^2\,e^2)$ in the secular terms %, the order of $\mathcal{O}(J_2\,e^2)$ in the long-period transformation equations,
and the $\mathcal{O}(J_2\,e^4)$ in the transformation equations. 
%The necessary equations are summarized in \ref{s:bgs}, where  the most relevant second-order corrections to the short-period terms are also provided, which will be eventually used for further comparisons.
\par

A set of numerical experiments that has been carried out to compare the relative efficiency of both approximate analytical solutions, Brouwer and DRI, is summarized in what follows. The reference orbits used as true solutions in the computation of the errors are obtained by the numerical integration of the main problem in Cartesian coordinates. In particular, the tests are based on a set of LEO with the following orbital elements \citep{GurfilLara2014}
\begin{equation} \label{iicc}
a=7000\,\mathrm{km}, \quad \omega=10^\circ, \quad \Omega=0, \quad f=15^\circ,
\end{equation}
and varying inclinations from equatorial to polar. The tests have been carried out for the case of almost-circular orbits, with $e=0.005$, and also, for  $e=0.075$, which is assumed to be the maximum eccentricity before re-entry for this kind of orbit.
\par

\begin{figure}[htbp]
\centerline{
\includegraphics[scale=0.79]{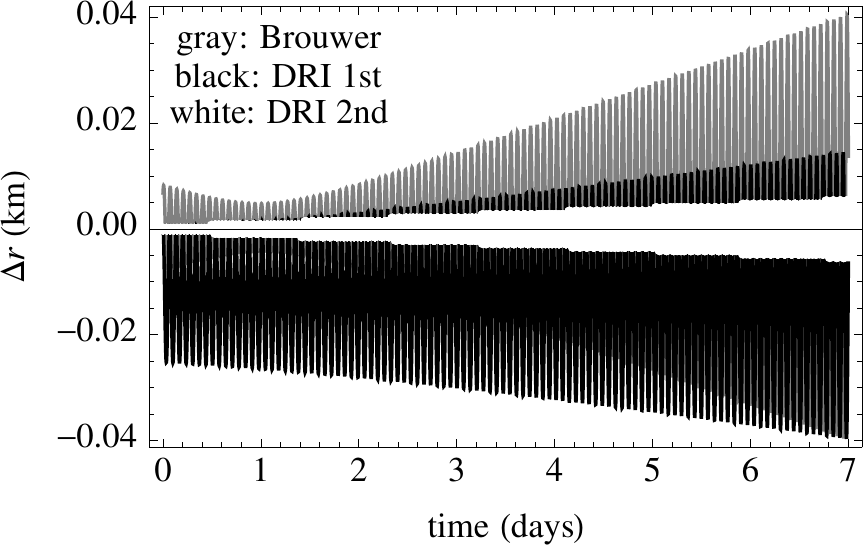} \quad \includegraphics[scale=0.75]{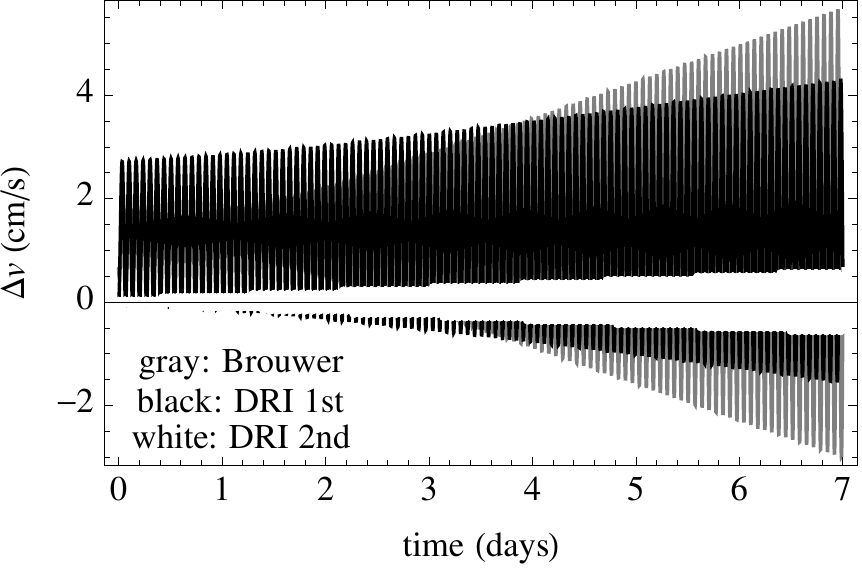} \\[1.5ex]
}
\centerline{
\includegraphics[scale=0.79]{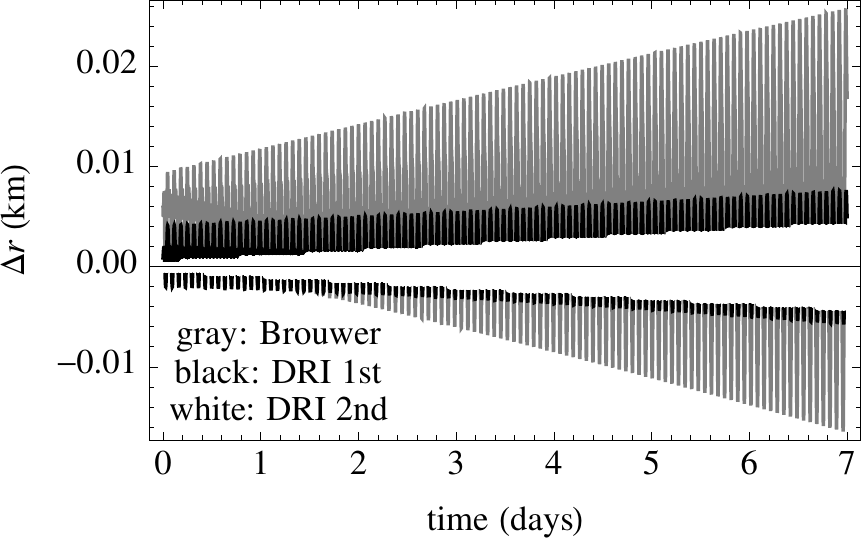} \quad \includegraphics[scale=0.75]{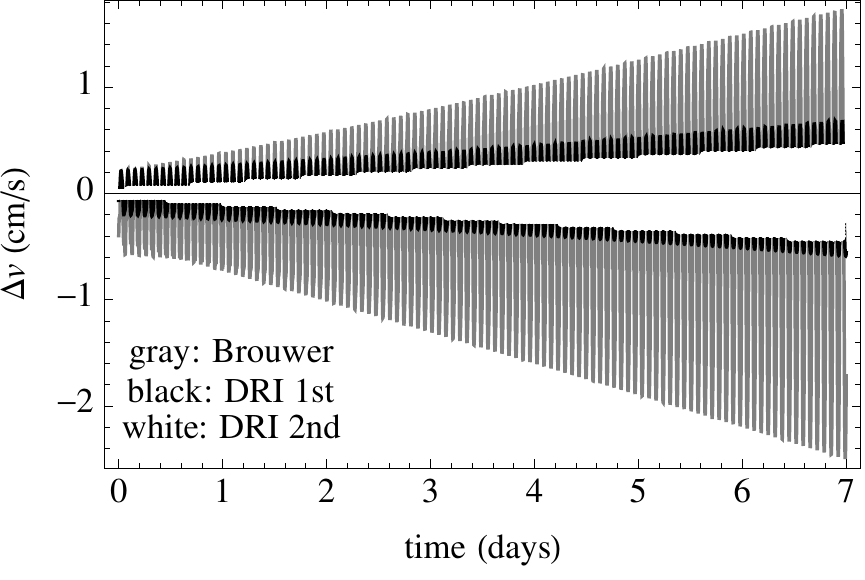} \\[1.5ex]
}
\centerline{
\includegraphics[scale=0.79]{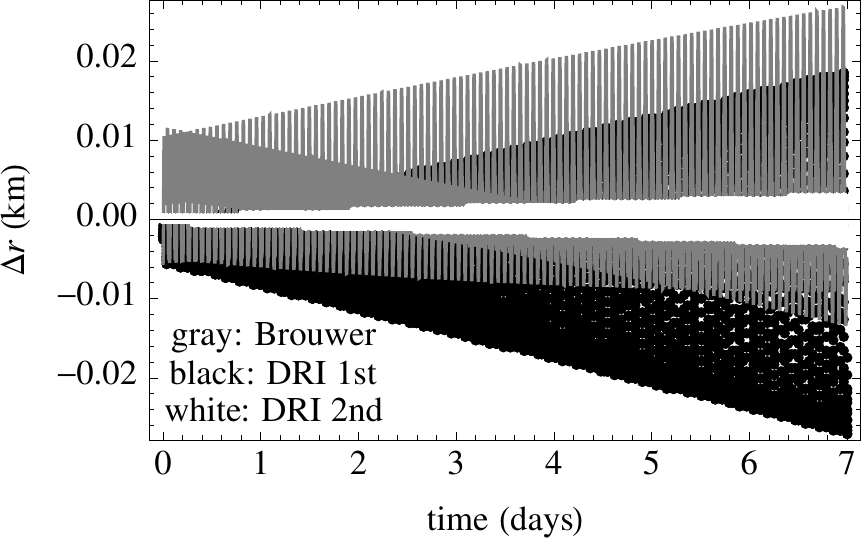} \quad \includegraphics[scale=0.75]{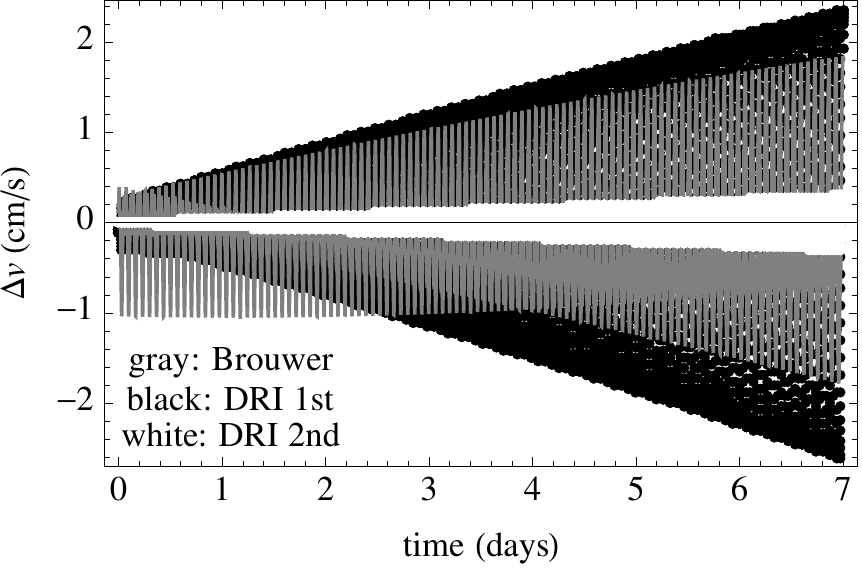} 
}
\caption{Propagation errors of Brouwer's solution (gray), first order of DRI (black) and second order of DRI (white, over the axis of abscissas). From top to bottom, $i=5$, $55$ and $89\deg$. Initial conditions in Eq.~(\ref{iicc}) for an almost-circular orbit with $e=0.005$}
\label{f:e005days7}
\end{figure}

\begin{figure}[htbp]
\centerline{
\includegraphics[scale=0.78]{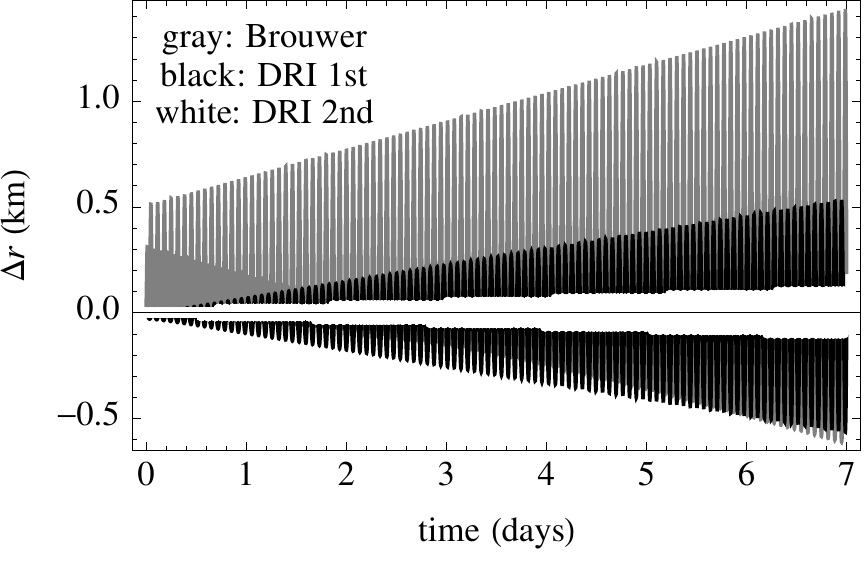} \quad \includegraphics[scale=0.77]{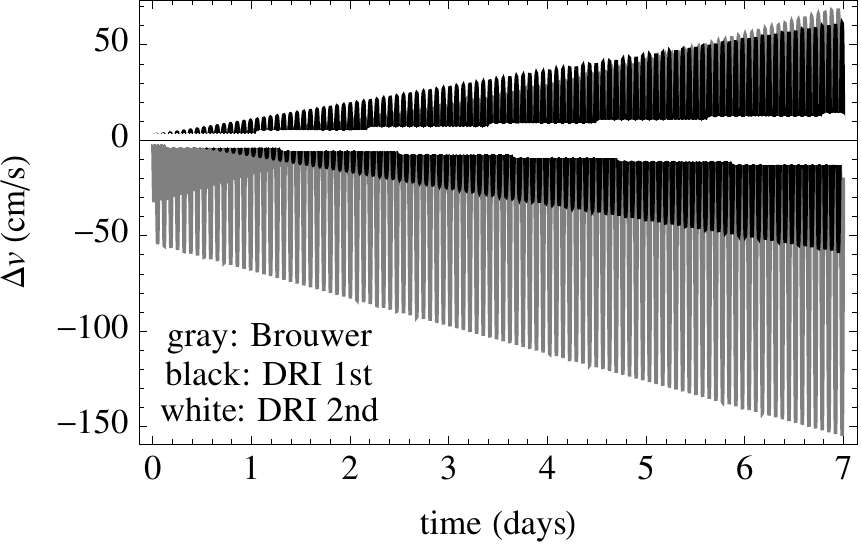} \\[1.5ex]
}
\centerline{
\includegraphics[scale=0.78]{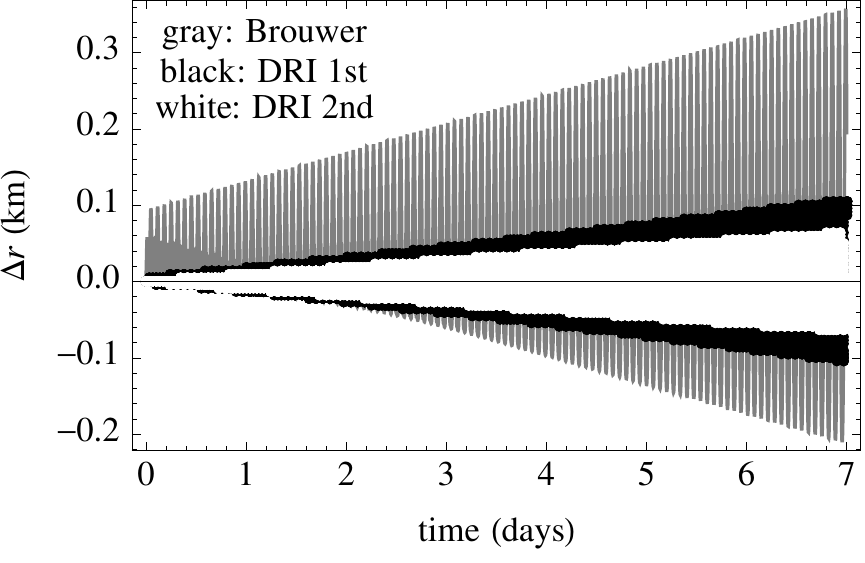} \quad \includegraphics[scale=0.77]{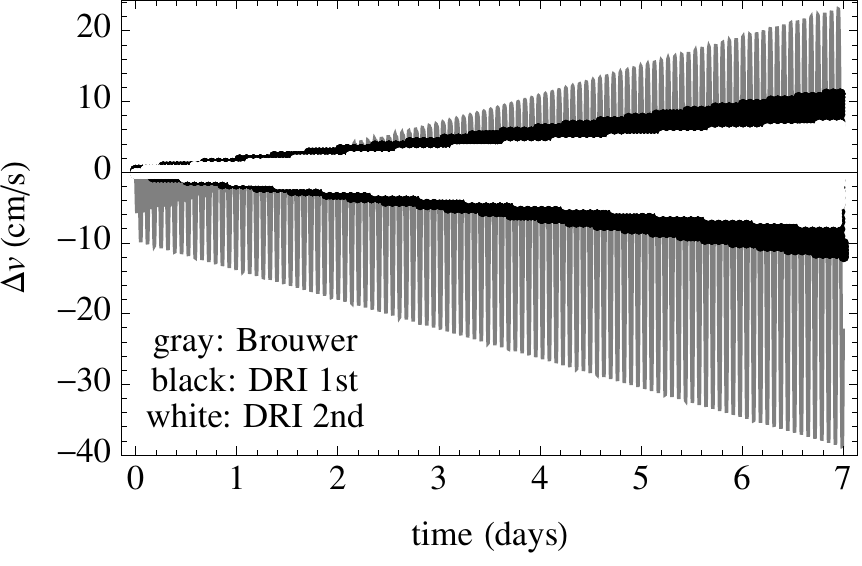} \\[1.5ex]
}
\centerline{
\includegraphics[scale=0.78]{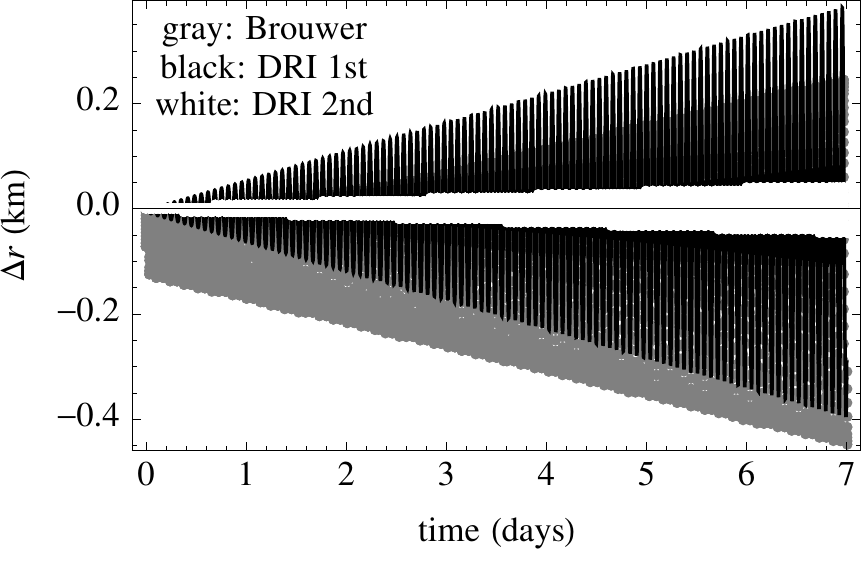} \quad \includegraphics[scale=0.77]{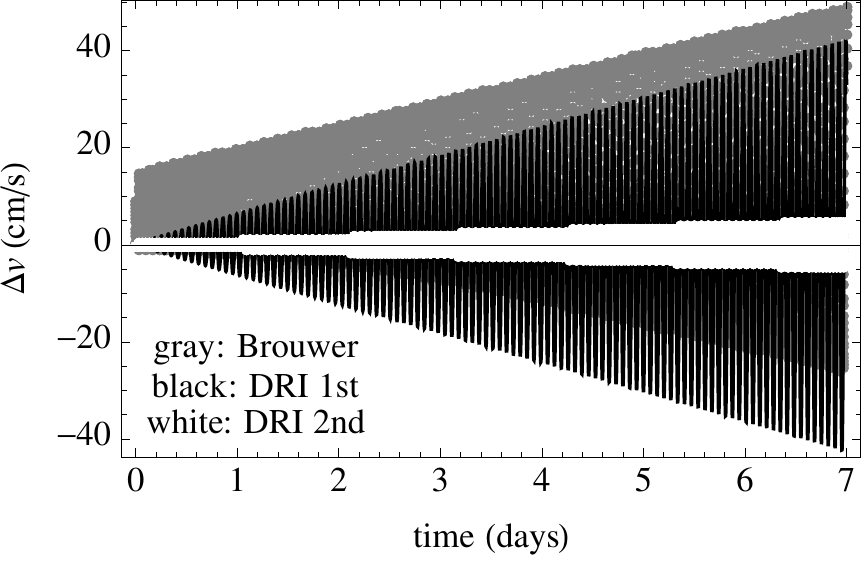} 
}
\caption{Propagation errors of Brouwer's solution (gray), first order of DRI (black) and second order of DRI (white, over the axis of abscissas), for, from top to bottom, $i=5$, $55$ and $89\deg$. Initial conditions in Eq.~(\ref{iicc}) for a slightly elliptic orbit with $e=0.075$}
\label{f:e075days7}
\end{figure}

In general, it is found that DRI reduces the errors by about one order of magnitude when compared with Brouwer's solution. Because the purpose of this research is to asses the relative efficiency of both different solutions, it is illustrative enough to compare the errors in the radial distance and in the modulus of the velocity vector. These comparisons are displayed in Figs.~\ref{f:e005days7} and \ref{f:e075days7}, which show the one-week time history of the errors obtained when using the different analytical solutions for orbits with initial conditions in Eq.~(\ref{iicc}). 
\par


Errors obtained when using the first order truncation of DRI, hereafter named DRI1, are also provided (black curves in Figs.~\ref{f:e005days7} and \ref{f:e075days7}). As observed in Fig.~\ref{f:e005days7}, in the case of the lower eccentricities DRI1 enjoys similar performance to Brouwer's solution for low- and high-inclination orbits. For medium inclination orbits, DRI1 clearly outperforms Brouwer's solution, and only performs slightly worse than DRI. The behavior of DRI1 remains the same for the higher-eccentricity orbits, as shown in Fig.~\ref{f:e075days7}, still performing better than Brouwer's solution in general, the case of high-inclination orbits being the unique where the performance of Brouwer's solution is comparable to that of DRI1.
\par

Similar computations have been carried out for time intervals spanning up to one month, which show analogous results. Thus, for 30 days propagation DRI always remain within less than 20 meter from the true (numerically computed) distance, and within less than 2\;cm/s of the true velocity for the lower-eccentricity orbits. Corresponding figures are 80\;m for distance and 4\;cm/s for velocity in the case of Brouwer's solution. These errors notably increase for higher eccentricities, but still errors obtained with DRI computations remain better than half a km in distance in all cases checked, while oscillations of the velocity errors remain below $\pm50$\;cm/s. In the case of Brouwer's solution, the errors after one month propagation may reach 5\;km in distance and 4\;m/s in velocity. Tests carried out include orbits with eccentricities in the range $0<e<0.1$.
\par

In addition to the higher precision obtained with DRI, because DRI equations are much simpler than Brouwer's equations, in all the cases tested the evaluation of DRI only requires about one fourth of the evaluation time required by Brouwer's solution, on average, in this way providing a definitive advantage of DRI over Brouwer's solution. 
\par

\begin{figure}[htbp]
\centerline{
\includegraphics[scale=1.075]{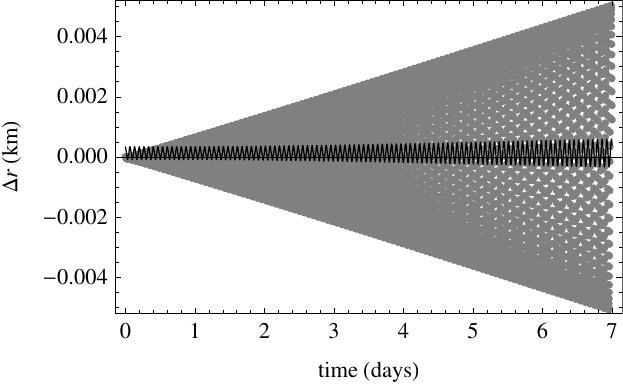} \quad \includegraphics[scale=1.07]{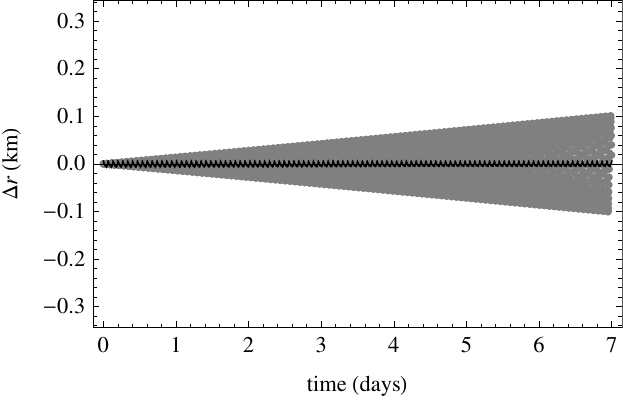} \\[1.5ex]
} \centerline{
\includegraphics[scale=1.075]{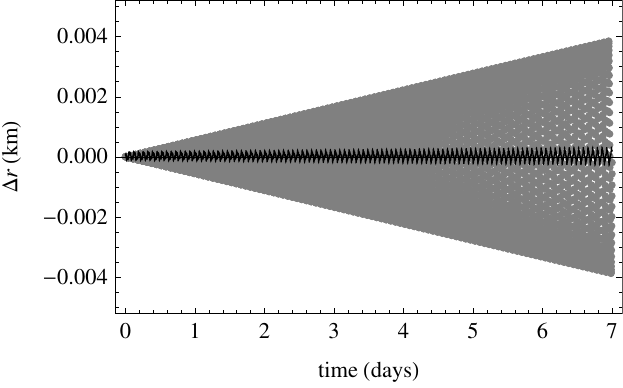} \quad \includegraphics[scale=1.07]{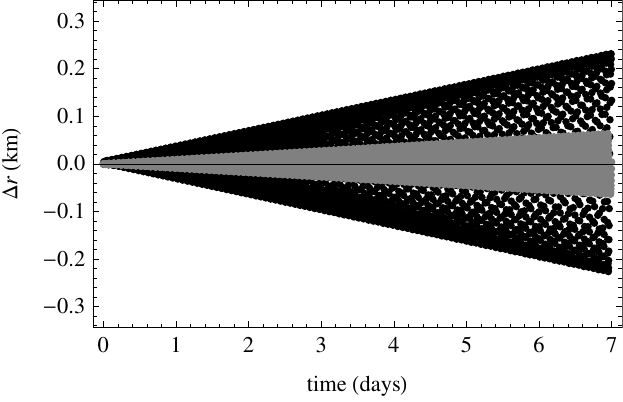} \\[1.5ex]
} \centerline{
\includegraphics[scale=1.075]{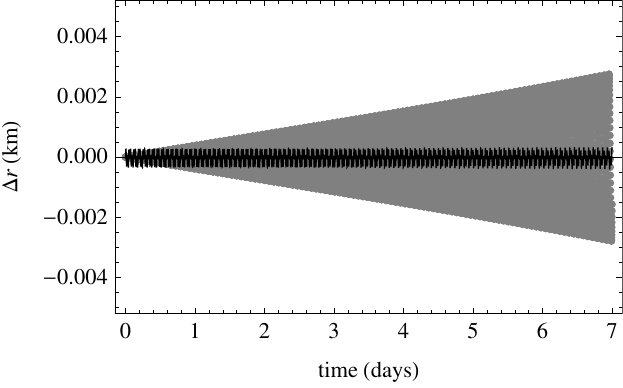} \quad \includegraphics[scale=1.07]{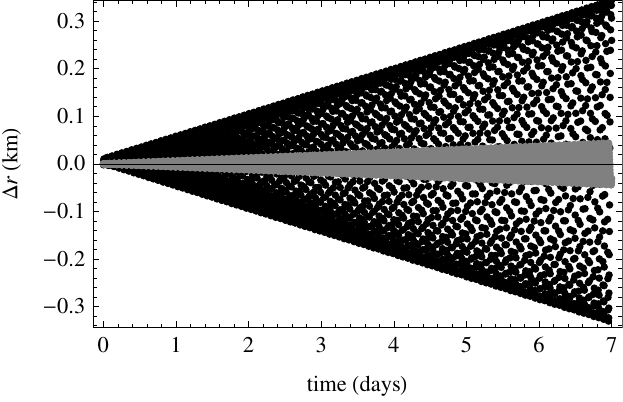} 
}
\caption{Propagation errors of the second order of DRI (gray) and Brouwer's solution with second order corrections for the short-period effects (black). From top to bottom, $i=5$, $55$ and $89\deg$. The left column corresponds to $e=0.005$ and the right one to $e=0.075$.}
\label{f:brouwer2nd}
\end{figure}

Results in Figs.~\ref{f:e005days7} and \ref{f:e075days7} clearly show the importance of taking into account second-order effects in the short-period corrections. To further illustrate this issue, a version of Brouwer's solution which includes the more relevant terms of the second-order corrections of the short-period terms has also been computed \citep[see][for instance]{DepritRom1970}. Comparisons of this improved Brouwer's solution with DRI are presented in Fig.~\ref{f:brouwer2nd} limited to errors in distance. Now Brouwer's improved solution clearly outperforms DRI for the lower eccentricity orbits, but this solution notably deteriorates for high-inclination orbits when the ellipticity is higher. Besides, Brouwer's improved solution unavoidably requires yet more computational effort than the original one.

\section{Conclusions}

When using Brouwer's gravity solution (in nonsingular elements) for short-term propagation, the more important source of errors comes from neglecting the contribution of second-order short-period effects in the computation of the theory, which in the case of LEO may result in the introduction of a concomitant uncertainty of tens of meters in the propagation of initial conditions. Besides, the computation of short-period corrections requires the evaluation of long Fourier series even when limited to first-order effects. Hence,  in order to speed computations, practical implementations of Brouwer's theory, as in SGP4, make additional simplifications by neglecting some of the first-order short-period corrections.
\par

On the contrary, first- and second-order short-period corrections are both obtained in a compact form of fast and straightforward evaluation when the gravity solution is approached with DRI, albeit the contribution of second-order secular and long-period terms is limited in DRI to the case of the lower eccentricity orbits. This fact makes that using DRI in the propagation of LEO leads to an increased precision when compared to the classical Brouwer's gravity solution, with the additional bonus of dramatically reducing the computation time. Hence, performances of existing orbit propagation software would benefit of replacing parts of the standard perturbation solution by analogous results based on intermediary orbits.
\par

While the model has been limited in this research to the main problem of artificial satellite theory, extension of the intermediary solution to other zonal models is encouraged by the quality of current results. The same assumptions made for integrating DRI as a quasi-Keplerian system apply also to models involving even zonal harmonics. In the case of odd zonal harmonics, separability of the intermediary obtained after the elimination of the parallax requires  making a stronger assumption by neglecting also long-period terms of the second order that are multiplied by the eccentricity. Nevertheless, the intermediary solution may still provide an acceptable accuracy for short-term propagation. The construction of a radial intermediary for a more complete zonal model is in progress, and corresponding results will be published elsewhere.

\appendix

\section{Solution of DRI} \label{a:DRI}

\subsection{Integration of the quasi-Keplerian system} \label{s:qK}
The quasi-Keplerian system is realized in the prime space; hence all variables of this appendix should be understood as prime variables.
\par
From initial conditions $(r_0,\theta_0,\nu_0,R_0,\Theta,N)$, where $\Theta$ and $N$ are constant in the prime space, the auxiliary constants
\begin{align} \allowdisplaybreaks
c \;=&\; \frac{N}{\Theta} \displaybreak[0] \\
\varepsilon \;=&\;-\frac{1}{4}J_2\,\frac{\alpha^2}{(\Theta^2/\mu)^2}, \displaybreak[0] \\
\tilde{\Theta} \;=&\; \Theta\,\sqrt{1-(2-6c^2)\,\varepsilon+(1-21c^4)\,\varepsilon^2}  \displaybreak[0] \\
\zeta \;=&\; \frac{\Theta}{\tilde\Theta}\left[1+\left(2-12c^2\right) \varepsilon-\left(3-105 c^4\right) \varepsilon ^2\right],
%=1+3\,(1-5c^2)\,\epsilon-3c^2\,(9-55c^2)\,\epsilon^2+\mathcal{O}(\varepsilon^3)
 \displaybreak[0] \\
\chi \;=&\; 6\varepsilon\left(1-7\varepsilon\,c^2\right)\frac{N}{\tilde\Theta},
%=6c\,\epsilon\left[1+\left(1-10c^2\right) \epsilon\right]+\mathcal{O}(\varepsilon^3),
 \displaybreak[0] \\
\tilde{p} \;=&\; \frac{\tilde{\Theta}^2}{\mu} \displaybreak[0] \\
a \;=&\; -\frac{\mu}{R_0^2+(\tilde{\Theta}^2/r_0^2) - 2\mu/r_0} \displaybreak[0] \\
e \;=&\; \sqrt{1-\frac{\tilde{p}}{a}}, 
\end{align}
as well as the initial values of the anomalies
\begin{align} \allowdisplaybreaks
f_0 \;\phantom{=}&\; \mathrm{from:} \qquad e\cos{f}_0=\frac{\tilde{p}}{r_0}-1, \quad e\sin{f}_0= R_0\,\sqrt{\frac{\tilde{p}}{\mu}} \displaybreak[0] \\
u_0 \;=&\; 2\arctan\left(\sqrt{\frac{1 - e}{1 + e}}\tan\frac{f_0}{2}\right), \displaybreak[0] \\
\ell_0 \;=&\; u_0-e\sin{u}_0,
\end{align}
are first computed.
Then, at a given time $t$, the anomalies
\begin{align} \allowdisplaybreaks
\ell \;=&\; \ell_0+ \sqrt{\frac{\mu}{a^3}}\,t, \displaybreak[0] \\
u \;\phantom{=}&\; \mathrm{from:} \qquad\ell=u-e\sin{u}, \displaybreak[0] \\
f \;=&\; 2\arctan\left(\sqrt{\frac{1+e}{1-e}}\tan\frac{u}{2}\right),
\end{align}
are computed first, and, then, the polar-nodal (prime) elements are obtained from the inverse sequence,
\begin{align} \allowdisplaybreaks
r(t) \;=&\; a\,(1-e\cos{u}), \displaybreak[0] \\
\theta(t) \;=&\; \theta_0+\zeta\,(f-f_0), \displaybreak[0] \\
\nu(t) \;=&\; \nu_0+\chi\,(f-f_0), \displaybreak[0] \\
R(t) \;=&\; \frac{\mu}{\tilde\Theta}\,e\sin{f},
%=\frac{\mu}{\Theta}\,e\sin{f}\frac{\Theta}{\tilde\Theta},
\end{align}
which is completed with the constants $\Theta$ and $N$.

\subsection{Short-period transformation} \label{s:shp}

The transformation from prime polar-nodal variables to original ones, and vice-versa, requires the computation of the total corrections
\begin{equation} \label{pncorrections}
\Delta_T\xi=\delta\,\Delta_1\xi+\frac{1}{2}\delta^2\,\Delta_2\xi,
\qquad
\delta=-\frac{1}{2}J_2\,\frac{\alpha^2}{p^2},
\end{equation}
where $\xi$ denotes any of the polar-nodal variables $\xi\in(r,\theta,\nu,R,\Theta,N)$, and $\delta=-\frac{1}{2}J_2\,(\alpha/p)^2$. Each total correction $\Delta_T\xi$ is to be added to the corresponding polar-nodal variable $\xi$. 
\par

The necessary corrections $\Delta_1\xi$ and $\Delta_2\xi$ are given below, where we recall that $\kappa\equiv\kappa(r,\Theta)$ and $\sigma\equiv\sigma(R,\Theta)$ are given in Eq.~(\ref{ef}). These corrections, as well as $\delta$, must be expressed in prime variables for the direct transformation $\Delta_{TD}\xi=\xi-\xi'$, and in original variables for the inverse transformation $\Delta_{TI}\xi=\xi'-\xi$.
\par

Remark that the first order corrections have been arranged in a form different from the original derivations of \citet{Deprit1981}, which is amenable for faster evaluation \citep{GurfilLara2014}. The second order corrections have been specifically derived and conveniently arranged for fast evaluation in the present work. For this sake, terms of the order of $e^2$ have been neglected from the expressions provided for the second order corrections, in agreement with the assumption made for the separability of DRI.

\subsubsection{First order corrections}

The first-order corrections are the same both for the direct and inverse transformations:
\begin{align} \allowdisplaybreaks
\Delta_1r \;=&\; p\left(1-\frac{3}{2}s^2-\frac{1}{2}s^2\cos2\theta\right), 
\displaybreak[0] \\[1ex]
\Delta_1\theta \;=&\; \left[\frac{3}{2}-\frac{7}{4}s^2+ (2-3s^2)\,\kappa\right]\sin2\theta \\ \nonumber
& -\left[5-6s^2+(1-2s^2)\cos2\theta\right]\sigma,
\displaybreak[0] \\[1ex]
\Delta_1\nu \;=&\;  c\left[(3+\cos2\theta)\,\sigma-\left(\frac{3}{2}+2\kappa\right) \sin2\theta\right],
\displaybreak[0] \\[1ex]
\Delta_1R \;=&\; \frac{\Theta}{p}\,(1+\kappa)^2s^2\sin2\theta,
\displaybreak[0] \\[1ex]
\Delta_1\Theta \;=&\; -\Theta\,s^2\left[\left(\frac{3}{2}+2\kappa\right)\cos2\theta+\sigma\sin2\theta\right],
\displaybreak[0] \\[1ex]
\Delta_1N \;=&\; 0,
\end{align}
where the right member must be expressed in prime variables for the direct transformation, and in original variables in the case of the inverse one.

\subsubsection{Second order corrections of the direct transformation}
Terms of $\mathcal{O}(e^2)$ have been neglected from the corrections. The right member must be expressed in prime variables.
\begin{align} \allowdisplaybreaks
\Delta_2r
\;= &\; \left\{-8+15 s^2-\frac{23}{4}s^4+\left(-\frac{3}{2}+\frac{7}{2}s^2-\frac{41}{16}s^4\right)\kappa \right. \\ \nonumber
& -\left[13-14s^2-\left(\frac{65}{8}-\frac{153}{16}s^2\right)\kappa\right]s^2\cos2\theta
-\left(\frac{1}{4}-\frac{1}{16}\kappa\right)s^4\cos4\theta \\ \nonumber
& \left. +\left[\left(\frac{27}{8}-\frac{51}{16}s^2\right)s^2\sin2\theta+\frac{9}{32}s^4\sin4\theta\right]\sigma \right\}p \displaybreak[0]
\\[1ex]
\Delta_2\theta
\;= &\; \left[8-29s^2+\frac{85}{4}s^4+\left(32-\frac{803}{4}s^2+\frac{1419}{8}s^4\right)\kappa\right]\sin2\theta \\ \nonumber
& +\left[\frac{9}{4}-\frac{3}{8}s^2-\frac{17}{8}s^4+\left(6-3s^2-\frac{55}{16}s^4\right)\kappa\right]\sin4\theta \\ \nonumber
&+  \left[72-121s^2+\frac{327}{8}s^4+ \left(-56+\frac{989}{4}s^2-\frac{1609}{8}s^4\right)\cos2\theta \right. \\ \nonumber
& \left. +\left(-3+3s^2+\frac{1}{8}s^4\right)\cos4\theta \right]\sigma 
\displaybreak[0] \\[1ex]
\Delta_2\nu \;= &\; \left\{ 
\left[\left(56-92s^2\right)\cos2\theta+\left(3-\frac{3}{2}s^2\right)(-9+\cos4\theta)\right]\sigma
   \right.\\ \nonumber
& \left. -\left[8-21s^2+\left(32-76s^2\right)\kappa\right]\sin2\theta
-\left(\frac{9}{4}+\frac{3}{4}s^2+6\kappa\right)\sin4\theta \right\}c 
\displaybreak[0] \\[1ex]
\Delta_2R \;= &\; 
\left\{ \left[16-16s^2+\left(\frac{237}{8}-\frac{437}{16}s^2\right)\kappa\right]s^2\sin2\theta \right. \\ \nonumber
& +\left(1+\frac{65}{32}\kappa\right)s^4\sin4\theta 
+\left[-\frac{3}{2}-\frac{1}{2}s^2+\frac{71}{16}s^4 \right. \\ \nonumber
& \left. \left. 
+\left(-\frac{95}{8}+\frac{231}{16}s^2\right) s^2 \cos2\theta
+\frac{17}{16}s^4\cos4\theta
\right]\sigma\right\}\frac{\Theta}{p} 
\displaybreak[0] \\[1ex]
\Delta_2\Theta \;= &\; \Theta  \left\{\left(\frac{9}{2}-\frac{25}{4}s^2+6\left(2-3s^2\right)\kappa\right)s^2 \right. \\ \nonumber
& -\left[8-\frac{15}{2}s^2+32\left(1-s^2\right)\kappa\right]s^2\cos2\theta
-\frac{3}{4}s^4\cos4\theta \\ \nonumber
& \left. +\sigma  \left[\left(-56+64s^2\right)s^2\sin2\theta+\frac{3}{2}s^4\sin4\theta\right]
\right\}
\displaybreak[0] \\[1ex]
\Delta_2N =& 0.
\end{align}

\subsubsection{Second order corrections of the inverse transformation}

Terms of $\mathcal{O}(e^2)$ have been neglected from the corrections. The right member must remain in original variables
\begin{align} \allowdisplaybreaks
\Delta_2r
\;= &\; \left\{ 8-12s^2+s^4+\left(\frac{3}{2}+\frac{1}{2}s^2-\frac{71}{16}s^4\right)\kappa \right. \\ \nonumber
& +\left[28-32s^2+\left(\frac{95}{8}-\frac{231}{16}s^2\right)\kappa\right]s^2\cos2\theta 
 -\left(1+\frac{17}{16}\kappa\right)s^4\cos4\theta \\ \nonumber
& \left. +\left[\left(-\frac{27}{8}+\frac{51}{16}s^2\right)s^2\sin2\theta-\frac{9}{32}s^4\sin4\theta\right]\sigma \right\}p \displaybreak[0] \\[1ex]
\Delta_2\theta
\;= &\; \frac{9}{4}-\frac{15}{8}s^2+2s^4+\left(6-3s^2-\frac{25}{16}s^4\right)\kappa +\left [-12+31s^2-\frac{73}{4}s^4 \right. \\ \nonumber
& \left. +\left(-40+\frac{819}{4}s^2-\frac{1371}{8}s^4\right)\kappa \right]\sin2\theta +\left[ -72+116s^2-\frac{243}{8}s^4 \right. \\  \nonumber
& \left.+ \left(26-\frac{1029}{4}s^2+\frac{1993}{8}s^4 \right)\cos2\theta + \left(-3+\frac{43}{8}s^4\right)\cos4\theta \right]\sigma \displaybreak[0] \\[1ex]
\Delta_2\nu \;= &\; \left\{ \left[12-21s^2+\left(40-76s^2\right)\kappa\right]\sin2\theta
 -\left(\frac{9}{4}-\frac{3}{4}s^2+6\kappa\right)\sin4\theta \right. \\ \nonumber
& \left. +\left[27-\frac{27}{2}s^2 + \left(-26+92s^2\right)\cos2\theta + \left(3+\frac{3}{2}s^2\right)\cos4\theta\right]\sigma \right\}c\displaybreak[0] \\[1ex]
\Delta_2R \;= &\; 
\left\{\left[-20+22s^2 - \left(\frac{333}{8}-\frac{725}{16}s^2\right)\kappa\right]s^2\sin2\theta \right. \\ \nonumber
& +\left(1+\frac{95}{32}\kappa\right)s^4\sin4\theta +\left[\frac{3}{2}-\frac{7}{2}s^2+\frac{41}{16}s^4  \right. \\ \nonumber
& \left. \left.
+ \left(-\frac{65}{8}+\frac{153}{16}s^2\right)s^2\cos2\theta - \frac{1}{16}s^4\cos4\theta \right]\sigma \right\}\frac{\Theta}{p} \displaybreak[0] \\[1ex]
\Delta_2\Theta \;= &\; \Theta\left\{ \left[\frac{9}{2}-\frac{25}{4}s^2+\left(12-18s^2\right)\kappa\right]s^2 \right. \\ \nonumber
& +\left[12-\frac{27}{2}s^2+\left(40-44s^2\right)\kappa\right]s^2\cos2\theta + \frac{3}{4}s^4\cos4\theta \\ \nonumber
& \left. +\left[\left(26-28s^2\right)s^2\sin2\theta - \left(\frac{3}{2}+\frac{9}{4}\kappa\right)s^4\sin4\theta \right]\sigma \right\}
\displaybreak[0] \\[1ex]
\Delta_2N \;= &\; 0.
\end{align}

\end{document}